\documentclass[11pt]{article}

\usepackage[letterpaper,width=135mm,height=203mm]{geometry}

\usepackage{graphicx}

\newcommand{\dd}[1]{\textbf{\textit{#1}}}

\newtheorem{theorem}{Theorem}
\newtheorem{lemma}{Lemma}

\newenvironment{Proof}{\textsc{Proof.}}{\hspace{8mm}\QEDmark\smallskip}
\usepackage{latexsym}
\newcommand{\QEDmark}{\mbox{$\Box$}}

\makeatletter
\long\def\@caption#1[#2]#3{\begingroup \@parboxrestore 
\if@minipage \@setminipage \fi \normalsize \sffamily \@makecaption {\csname fnum@#1\endcsname }{\ignorespaces #3}\par \endgroup}
\makeatother

\newcommand{\upp}{{\mathit{c}}}
\newcommand{\prop}{{\mathit{cp}}}

\begin{document}

\begin{center} 
{\bfseries \Large Bounds on Coloring Trees without Rainbow Paths} \\[2mm]
Wayne Goddard, Tyler Herrman, Simon J. Hughes \\[1mm]
School of Computing \& School of Mathematical and Statistical Sciences \\
Clemson University, Clemson SC USA
\end{center}

\begin{abstract}
For a graph with colored vertices, 
a rainbow subgraph is one where all vertices 
have different colors.
For graph $G$, let $\upp_k(G)$ denote the maximum number of different colors in a coloring 
without a rainbow path on $k$ vertices, and 
 $\prop_k(G)$ the maximum number of colors if the coloring is required to be proper. 
 The parameter $\upp_3$ has been studied by multiple authors. 
 We investigate these parameters for trees and $k \ge 4$. We first calculate them
 when $G$ is a path, and determine when the optimal coloring is unique.
 Then for trees $T$ of order $n$, we show that the minimum value of $\upp_4(T)$ and $\prop_4(T)$ is $(n+2)/2$, and the trees with the minimum value 
of $\prop_4(T)$ are  the coronas. Further, 
the minimum value of $\upp_5(T)$ and $\prop_5(T)$ is $(n+3)/2$ ,
and the trees with the minimum value of either parameter are octopuses.
\end{abstract}

\section{Introduction}

We consider undirected graphs where the vertices receive colors. 
We define a subgraph as \dd{rainbow} if all its vertices 
receive different colors, and we study colorings where for some fixed graph $H$ there
is no rainbow subgraph isomorphic to $H$. This question was first studied
for the path on three vertices by Bujt\'as et al.~\cite{BSTSD} and for stars in general by Bujt\'as et al.~\cite{BSTDP},
and then studied in \cite{CDHH,GXfirst,GXsecond} inter alia. There has also been work 
on the case where $H$ must be induced~\cite{AM}. The problem is also a  special case of more general
questions introduced and studied earlier by Voloshin~\cite{Voloshin}. (The edge-coloring
version is much more studied, where it is called anti-Ramsey theory.)

Our focus here is on the case that $H$ is a path. And specifically on the maximum number of
colors one can use on a graph and there not be a rainbow path.
For graph $G$, let $\upp_k(G)$ denote the maximum number of different colors one can 
use without there being a rainbow $P_k$ (meaning a path with $k$ vertices), where
the coloring is not required to be proper. Let $\prop_k(G)$ denote the maximum number
of colors with the additional constraint that adjacent vertices receive different colors; that is, it is a proper coloring.
Note that $\prop_k(G)$ might not exist; for example, it does not exist for the complete graph $K_n$ where $n\ge k$.

As noted above, the function $\upp_3(G)$ has already been studied. See for example \cite{BSTSD,CDHH,GXfirst} and the references therein. 
The equivalent $\prop_3(G)$ is uninteresting: the only way a $P_3$ can be properly colored 
without being rainbow is that the first and third vertex have the same color; so such a coloring of $G$ exists
only when $G$ is bipartite. The parameter $\upp_4(G)$ is also briefly studied in~\cite{GXsecond}.

We proceed as follows. In Section~\ref{s:paths} we consider the colorings of paths and determine when the extremal colorings are unique.
In Section~\ref{s:p4} we show that the minimum value of $\upp_4(T)$ and $\prop_4(T)$  for trees of  order $n$ is $(n+2)/2$, and the trees with the minimum value 
of $\prop_4(T)$ are precisely the coronas.
In Section~\ref{s:p5} we show that the minimum value of $\upp_5(T)$ and $\prop_5(T)$ for trees of  order $n$ is $(n+3)/2$, 
and the trees with the minimum value of either parameter are octopuses.
In Section~\ref{s:end} we conclude with brief thoughts for future research

\section{Colorings Paths Without Rainbow Paths}  \label{s:paths}

We begin with the calculation of the parameters for paths and determining when the optimal coloring is unique. 
There has been previous work. Observation 19 of~\cite{GWX}
considered the problem of coloring the path such that there is no specified rainbow subpath and no three 
consecutive vertices
receive the same color; since the optimal coloring without a rainbow path clearly does not
have three consecutive vertices of the same color, the formula given there applies to $\upp_k(P_n)$. The same formula can
also be read out of the results on mixed interval hypergraphs in~\cite{BV}. For proper colorings,  
the special case of $\prop_4(P_n)$ was resolved in~\cite{GXsecond}.
Nevertheless, we include proofs of all formulas, since this proof enables us to determine when the extremal coloring is unique.

For a graph $G$ and vertex $w$, we define the operation of \dd{attaching a $P_m$} as adding a copy of $P_m$ and 
joining one end of the $P_m$  to~$w$. See Figure~\ref{f:attach} for an example. We will need the following lemma, which is possibly interesting in its own right:

\begin{figure}[h]
\centerline{\includegraphics{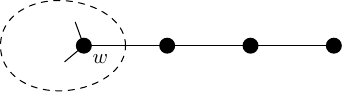}}
\caption{Attaching a $P_3$}
\label{f:attach}
\end{figure}

\begin{lemma} \label{l:addition}
(a) Assume $k\ge 2$. For any graph $G$ and vertex $w$, if graph $G_1$ is obtained from $G$ by attaching $P_{k-1}$ 
to~$w$, then $\upp_{k} ( G_1 ) = \upp_k (G ) +  k- 2$. \\[1mm]
(b) Assume $k\ge 3$. For any graph $G$ and \dd{end-vertex} $w$, if graph $G_2$ is obtained from $G$ by attaching $P_{k-2}$ 
to $w$, then $\prop_{k} ( G_2 ) = \prop_k (G ) +  k- 3$. 
\end{lemma}
\begin{Proof}
(a) Let $X$ denote the attached $P_{k-1}$.  It cannot happen that
every vertex in $X$ gets its own unique color, since that creates a rainbow $P_k$ with $w$.
On the other hand, any coloring of~$G$ is extendable to $G_1$ 
by giving the first vertex of $X$ the same color as $w$ 
and giving the remaining $k-2$ vertices of $X$ each their own unique color.
This proves the formula.

(b)  Let $Y$ denote the attached $P_{k-2}$.   It cannot happen that
every vertex in~$Y$ gets its own unique color, since that creates a rainbow $P_k$ with $w$ and $w$'s neighbor
in $G$ (which necessarily have different colors). On the other hand, any coloring of $G$ is extendable to $G_2$ 
by giving the first vertex of $Y$ the same color as $w$'s neighbor in $G$, 
and giving the remaining $k-3$ vertices of $Y$ each their own unique color.
This proves the formula.
\end{Proof}

Note that part (b) of Lemma~\ref{l:addition} does not generalize to all $w$: for example, if $G=K_3$ then $\prop_4(G)=3$ but $\prop_4(G_2)=3$ too.
(It is of course true that $\prop_{k} ( G_2 ) \le \prop_k (G ) +  k- 3$ for all $w$.) Lemma~\ref{l:addition} enables us to prove the following:

\begin{theorem}  \label{t:paths} Let $k\ge 2$. \\
(a) For $n\ge 1$ it holds that $\upp_k( P_n ) = \lfloor (k-2)n/(k-1) \rfloor +1$.  \\
(b) For $n\ge k-1$ the optimal coloring is unique exactly when $n$ is a multiple of $k-1$.
\end{theorem}
\begin{Proof}
(a) For fixed $k$ we prove the bound by induction on $n$. 
For the base case, note that the formula gives $n$ for $n\le k-1$, which is correct. For $n\ge k$,
the above lemma implies that $\upp_k( P_n ) = \upp_k( P_{n-k+1} ) + k-2$, and so
the formula follows from the inductive hypothesis. For example, 
 for $P_{11}$ without rainbow $P_5$, an optimal coloring uses nine colors such as 
\begin{quote}
\verb|12344567789| .
\end{quote}

(b) Uniqueness also follows by induction. The base case is the range $k-1 \le n\le 2k-3$.
For $n=k-1$ the optimal coloring is rainbow. Otherwise
the optimal coloring has all vertices but two having their own color. Such a coloring
is valid if and only if 
the two vertices with the same color are contained in positions $n-k+1$ through $k$,
and thus the coloring is not unique. So assume $n\ge 2k-2$. 

To have equality in the bound, by the recurrence it follows that
 the coloring of $P_{n-k+1}$ must be optimal. On the other hand, any coloring of $P_{n-k+1}$ is extendable to $P_n$ 
by duplicating the color of the last vertex of $P_{n-k+1}$ and then giving the remaining $k-2$ vertices of $X$ each their own color.
Hence for the optimal coloring of $P_n$ to be unique, so must the optimal coloring of $P_{n-k+1}$ and thus by the induction hypothesis the divisibility
condition is necessary. To show that the condition is also sufficient, note that when the coloring of $P_{n-k+1}$ is unique, 
its last $k-1$ vertices have different colors. So the only way to get $k-2$ colors on $X$ 
is as described before: the first vertex of $X$ must have the same color as the last vertex of the $P_{n-k+1}$.
It follows that the optimal coloring of $P_n$ is unique.
\end{Proof}

\begin{theorem} \label{t:pathsProper} Let $k\ge 3$. \\
(a) For $n\ge 2$ it holds that $\prop_k( P_n ) =  \lfloor ((k-3)n+1)/(k-2) \rfloor +1$.  \\
(b) For $n\ge k-1$ the optimal coloring is unique exactly when $n$ is $1$ more than a multiple of $k-2$.
\end{theorem}
\begin{Proof}
(a) For fixed $k$ we prove the bound by induction on $n$. 
For the base case, note that the formula gives $n$ for $n\le k-1$, which is correct. For $n\ge k$,
the above lemma implies that $\prop_k( P_n ) = \prop_k( P_{n-k+2} ) + k-3$, and so
the formula follows from the induction hypothesis.
For example, for $P_{11}$ without rainbow $P_5$, the optimal proper coloring uses eight colors such as
\begin{quote}
\verb|12343565787| .
\end{quote}

(b) Uniqueness also follows by induction.
The base case is the range $k -1 \le n\le 2k-4$.
For $n=k-1$ the optimal coloring is rainbow. Otherwise
the optimal coloring has all vertices but two having their own color. Such a coloring
is valid if and only if 
the two vertices with the same color are contained in positions $n-k+1$ through $k$,
and thus the coloring is not unique. So assume $n\ge 2k-3$. 

To have equality in the bound, by the recurrence it follows that
the coloring of $P_{n-k+2}$ must be optimal. On the other hand, any coloring of $P_{n-k+2}$ is extendable to $P_n$ 
by duplicating the color of the penultimate vertex and then giving the remaining $k-3$ vertices of $Y$ each their own  color.
Hence for the coloring of $P_n$ to be unique, so must the coloring of $P_{n-k+2}$, and so by the induction hypothesis the divisibility
condition is necessary. To show that the condition is also sufficient, note that when the coloring of $P_{n-k+2}$ is unique, 
its last $k-1$ vertices have different colors. So the only way to get $k-3$ colors on~$Y$ 
is as described before: the first vertex of~$Y$ must have the same color as the penultimate vertex of the $P_{n-k+2}$.
It follows that the optimal coloring of $P_n$ is unique.
\end{Proof}


One can also derive the formula for $\upp_k(P_n)$ using the fact that, in trees, the  parameter
 is intimately related to the minimum number of edges one must remove to destroy all copies of $P_k$.
For $H$ a fixed graph, define a set $F$ of edges in a graph~$G$ as \dd{$H$-thwarting} if removing all of $F$ from the graph $G$
destroys all copies of~$H$. The $H$-thwarting number, $\theta_H(G)$,  is the minimum 
number of edges whose removal destroys all copies of $H$. 
In a coloring, we call an edge \dd{monochromatic} if its two ends have the same color.
Note that if the monochromatic edges form a $H$-thwarting set, then every $H$ contains
a monochromatic edge and hence the coloring is valid,
that is, has no rainbow $H$. Hence in general graphs $G$ there is the inequality $\upp_H(G) \ge n - \theta_H(G)$,
where $\upp_H(G)$ denotes the maximum number of colors in a coloring of $G$ without a rainbow copy of $H$.
Theorem~16 in the paper~\cite{GXfirst} showed that in trees there is equality:

\begin{theorem} \cite{GXfirst}
In any tree $T$ of order $n$, it holds that $\upp_H(T) = n-\theta_H(T)$.
\end{theorem}

\section{Coloring Trees Without a Rainbow $P_4$}  \label{s:p4}

A general lower bound for bipartite graphs in the $P_4$ case was obtained in  \cite{GXsecond}. Namely, it was 
observed that in a bipartite graph with bipartition $(X,Y)$,
if one gives each vertex in $X$
its own unique color and gives all the vertices in $Y$ the same color, 
the result is a valid proper coloring: every copy of $P_4$ has two vertices 
from $Y$. Since the bigger partite set has at least half the vertices, it follows that:

\begin{theorem} \cite{GXsecond} \label{t:first}
For any connected bipartite graph $G$ on $n\ge 2$ vertices it holds that 
\[
    \upp_4(G) \ge \prop_4(G) \ge \lceil n/2 \rceil +1 .
 \] 
 \end{theorem}
 
 It was also noted in \cite{GXsecond} that, if a graph $G$ has a perfect matching, then $\prop_4(G) \le n/2+1$.
 The value for paths given in Theorem~\ref{t:pathsProper} is thus recovered.
 
 \subsection{Trees with Extremal $\upp_4$}

\begin{figure}[h]
\centerline{\includegraphics{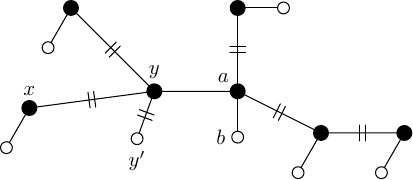}}
\caption{A corona and a minimum $P_4$-thwarting set}
\label{f:corona}
\end{figure}

Perhaps surprisingly, the paths do not have the minimum value of $\upp_4(T)$ for a given order.
Indeed we show that the trees $T$ with the minimum value of $\upp_4(T)$ are precisely the coronas.
The \dd{corona} of a graph is defined by taking the graph, and for each vertex~$w$
adding one new vertex, called a \dd{foot}, adjacent only to~$w$. (This doubles the number of vertices.)
The original graph is called the \dd{core}. 
Figure~\ref{f:corona}  gives an example where the core is a tree (and the notched edges form a minimum 
$P_4$-thwarting set).

We will need the following lemma:

\begin{lemma} \label{l:thwartLeaf}
If $T$ is a tree with end-vertex $y'$, then there exists a minimum $P_4$-thwarting set
that does not contain the edge incident with $y'$.
\end{lemma}
\begin{Proof}
Let $F$ be a minimum thwarting set.
Say $y'$ has neighbor $y$.
Assume the edge $yy'$ is in $F$. The minimality of $F$ 
means that if we put $yy'$ back into the graph there will be a $P_4$, say $y'yab$, where only edge $y'y$ is in $F$. 
It follows that every other edge incident with $y$ must be in $F$, else we have a $P_4$ starting $bay$.
So one can change $F$ by replacing $y'y$ by $ya$ and still have a minimum thwarting set of $T$.
\end{Proof}

\begin{theorem} \label{t:coronaGood}
For a corona $H$ derived from core tree $B$, it holds that
$\upp_4(H) = n/2 +1$, where $n$ is the order of $H$.
\end{theorem}
\begin{Proof}
It suffices to show that $\theta_{P_4} (H) = n/2-1$. Such a thwarting set
is achieved by taking all the edges in the core.
It remains to show that there is no smaller thwarting set.

We prove that $\theta_{P_4} \ge n/2-1$ by induction on the number of vertices of the core graph.
The base case of the induction is that $B$ has one vertex. That
corresponds to $H$ being $K_2$; this has $\theta_{P_4} =0$.
So assume $B$ has at least two vertices. 

Consider a vertex $x$ that is an end-vertex in $B$. Say its neighbor in $B$ is $y$.
By Lemma~\ref{l:thwartLeaf}, there is a minimum thwarting set $F$ that does not
contain the edge joining $y$ to its foot, say $y'$. It follows that some edge incident with $x$ is 
in~$F$. Let $T'$ be obtained from $T$ by removing $x$ and its foot neighbor.
Then the portion of~$F$ restricted to $T'$ is a thwarting set of $T'$, and hence by the induction hypothesis,
has size at least $(n-2)/2-1$. But $F$ also contains an edge incident with $x$, and so has
size at least $n/2-1$, as required.
\end{Proof}

We show next that coronas are the only examples of trees where $\upp_4 = n/2+1$.
We will need the following lemma.

\begin{lemma} \label{l:coronaThwart}
If $T$ is a corona, then there is a minimum $P_4$-thwarting set that contains
any one designated leaf edge.
\end{lemma}
\begin{Proof}
Assume the designated edge joins end-vertex $y'$ with neighbor $y$. Then a thwarting set can be
constructed by taking $yy'$ together with all edges of the core, except for one incident with $y$. (See Figure~\ref{f:corona} earlier.)
\end{Proof}

\begin{theorem} \label{t:notCorona}
If $T$ is a tree on $n\ge 2$ vertices that is not a corona, then $\upp_4(T) > n/2+1$.
\end{theorem}
\begin{Proof}
By Theorem~\ref{t:first},
we already know this for trees of odd order. So assume that $n$ is even.
It suffices to show that $\theta_{P_4}(T) < n/2-1$.
The proof is by induction on the diameter of $T$. 
If the diameter is $1$ then $n=2$ but the tree is a corona.
If the diameter is $2$, then the tree is a star and $\theta_4=0$; so the result is true.

So assume the diameter of the non-corona tree $T$ is at least $3$. Then since $P_4$ 
is  a corona, we know $T$ is not $P_4$ and hence $n\ge 6$.
Consider a longest (diametrical) path $Q$ in 
the tree. Say the path starts with vertices $abcd$.\smallskip 

\dd{Case 1:} \textit{Vertex $b$ has degree $2$.}
Then $T' = T-\{a,b\}$ is a tree. 
If $T'$ is not a corona, then $\theta_{P_4} (T') < (n-2)/2-1$ by the induction hypothesis, and we can 
extend to a thwarting set of $T$ by adding the edge $bc$.

So assume $T'$ is a corona. Then since $T$ is not 
a corona, it must be that $c$ is an end-vertex of $T'$. 
At the same time, by Lemma~\ref{l:coronaThwart}, there is a minimum thwarting set of $T'$ (that is, of size $(n-2)/2-1$)
that uses edge $cd$. This is also a thwarting set of $T$.\smallskip

\dd{Case 2:} \textit{Vertex $b$ has degree $r > 2$.}
Then by the maximality of the path~$Q$, 
every neighbor of $b$ except $c$ is an end-vertex. 
Let $T'$ be the graph obtained from~$T$ by deleting $b$ and all its end-vertex neighbors.
Note that $T'$ is a tree. By Theorem~\ref{t:first}, it has a thwarting set of size at most
$(n-r)/2-1$. 
One can extend that set to a thwarting set of $T$ of size $n/2-r/2 < n/2-1$ by adding edge $bc$.
\end{Proof}

At the other extreme, the question of the maximum value of $\upp_4$ is trivial, 
since the value is $n$ if and only of the graph does not contain a copy of $P_4$. 
One can also readily determine the trees where the value is $n-1$. 
We define a \dd{multi-corona} as a graph that results from adding one or more feet to every vertex of a base graph.

\begin{theorem} \label{t:upperN-1}
For a tree $T$ of order $n$ containing $P_4$, it holds that
$\upp_4(T) = n -1 $ if and only if $T$ is a subgraph of a multi-corona of $P_4$. 
\end{theorem}
\begin{Proof}
It is immediate that the middle edge of the core $P_4$ forms a thwarting set. 
On the other hand, consider a tree $T$ with a thwarting set consisting of only the edge $ab$.
Then all copies of $P_4$ in $T$ contain the edge $ab$. Thus vertex $a$ cannot have
two neighbors of degree more than $1$, and if it does have a neighbor $a'$ of degree more than $1$, that 
vertex can only have end-vertex neighbors. Similarly $b$ has at most one neighbor $b'$ of degree more than $1$.
Thus $T$ is a subgraph of a multi-corona of $a'abb'$.
\end{Proof}

 \subsection{Trees with Minimum $\prop_4$}

\begin{figure}[h]
\centerline{\includegraphics{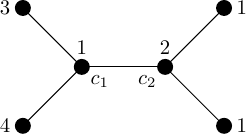}}
\caption{The double star $D_2$  with $\prop_4 = 4$}
\label{f:doubleStar}
\end{figure}

In contrast to the case for $\upp_4$, it seems that there is no simple description of the trees where $\prop_4 =  n/2 + 1$. We have seen that this is true of paths.
It is also true for the 
\dd{double star} $D_{b}$, defined by taking two stars each with $b$ end-vertices and
adding one edge joining the two centers $c_1$ and $c_2$. See Figure~\ref{f:doubleStar} for an example.

\begin{lemma}
For a double star $D_b$ it holds that $\prop_4 (D_b) = b + 2$.
\end{lemma}
\begin{Proof}
The lower bound follows from Theorem~\ref{t:first}.
Consider a valid coloring of~$D_b$. Then the central vertices, say $c_1$ and $c_2$, receive different 
colors. Consider some vertex with a third color;  say vertex $v$
adjacent to $c_1$. Then by the $vc_1c_2w$ path, every vertex $w$ adjacent to $c_2$
must have either the color of $v$ or the color of~$c_1$; that is, not a new color.
So the number of colors other than that of $c_1$ and $c_2$ is at most
the number of leaf-neighbors of~$c_1$, which is $b$.
Hence the double-star has $\prop_4 = b + 2$.
\end{Proof}

We next show that, if $T$ is a tree with $\prop_4(T) =  |T|/2 + 1$, then so is one with $K_2$ attached.
We observed earlier that part (b) of Lemma~\ref{l:addition} does not extend to general attachers $w$. 
However, it turns out that it does extend if the underlying graph is a tree, at least in the case $k=4$.

We will need the following idea. For a coloring,
define a vertex $x$ as \dd{boring} if either (i) all neighbors of $x$ have the
same color, or (ii) all vertices at distance~$2$ from $x$ have the same color as $x$, or
both. For example, in the coloring of $D_2$ in Figure~\ref{f:doubleStar} every vertex is boring.
We claim that in 
a tree $T$, one can choose an optimal $\prop_4$-coloring such that every vertex is
boring.

\begin{lemma} \label{l:boring}
If $T$ is a tree, then there exists an optimal $\prop_4$-coloring such that every vertex is
boring.
\end{lemma}
\begin{Proof}
Consider the optimal coloring of $T$ with the most boring vertices, and suppose
there is a vertex $x_3$ that is not boring. Then there is a vertex $x_1$ at
distance~$2$ from $x_3$ with a different color. Let $x_2$ be their common neighbor.
Say $x_1$ has color $1$, $x_2$ has color $2$, and $x_3$ has color $3$.
Further, there must be a neighbor of $x_3$, say $x_4$, that 
has a color different from $x_2$. 
Since $x_1x_2x_3x_4$ is not rainbow, vertex $x_4$ must have color~$1$.
Indeed, all neighbors of $x_3$ must have color $1$ or $2$.
Furthermore, all vertices at distance two must have color $1$, $2$, or~$3$.
See Figure~\ref{f:boringPic} for an example.

\begin{figure}[h]
\centerline{\includegraphics{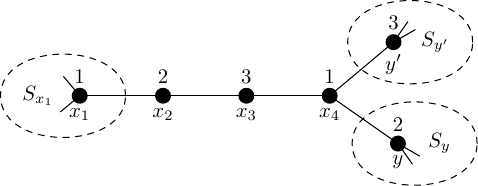}}
\caption{A possible coloring}
\label{f:boringPic}
\end{figure}

Now, we re-color the tree $T$. For each vertex $y$ at distance $2$ from $x_3$,
let $S_y$ be the subtree consisting of $y$ and all vertices whose path to $x_3$ 
goes via $y$. If $y$ has color $1$, then in $S_y$ change every vertex
with color $1$ to color $3$, and vice versa. 
If $y$ has color $2$, then in $S_y$ change every vertex
with color $2$ to color~$3$, and vice versa.
One does not lose any color in the process, as $x_2x_3x_4$ still contains
all three colors. 

We claim that the new coloring is still a valid coloring:
within $S_y$ only names of colors have been changed, and similarly
with including the common neighbor of $y$ and $x_3$. Further, any $P_4$ containing both
$y$ and~$x_3$, or containing both $y$ and another vertex at distance two from $x_3$,
has two vertices of color $3$. Finally, we note that $x_3$ is now boring (since all vertices at distance~$2$ have the same color),
as are all its neighbors. Further, every other vertex that was boring remains so.
This re-coloring increases the number of boring vertices, and so contradicts the choice of 
coloring. That is, the supposition that the coloring had a vertex that is not boring is false.
\end{Proof}

Using Lemma~\ref{l:boring} we can prove a result about attachments:

\begin{lemma} \label{l:treeAddition}
For any nontrivial tree $T$ and vertex $w$, let tree $T_2$ be obtained from $T$ by 
attaching $P_{2}$. Then $\prop_{4} ( T_2 ) = \prop_4 (T ) +  1$. 
\end{lemma}
\begin{Proof}
We noted earlier that $\prop_{4} ( T_2 ) \le \prop_4 (T ) +  1$, since the two new vertices cannot both get a unique color. So it remains
to find a suitable coloring.

By Lemma~\ref{l:boring} there exists an optimal coloring of $T$ where $w$ is boring.
We color the $P_2$ as follows.
If all neighbors of $w$ have the same color, say red, then 
color the first vertex of $P_2$ red and give the other vertex a new color.
If all vertices at distance $2$ from~$w$ have the same color as $w$, say blue,
then give the first vertex of $P_2$ a new color and color the 
other vertex blue. In either case the result is a valid coloring.
\end{Proof}

\begin{figure}[h]
\centerline{\includegraphics{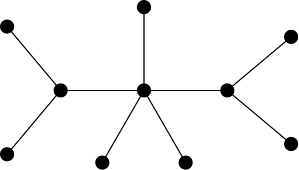}}
\caption{A tree with $\prop_4(T) = 6$}
\label{f:properSmall}
\end{figure}

But there are many other trees $T$ with $\prop_4(T) = |T|/2 + 1$. Figure~\ref{f:properSmall} shows an example.

At the other extreme, the question of the maximum value of $\prop_4$ can be also be considered. 
Using Theorem~\ref{t:upperN-1} one can show:

\begin{theorem}
For a tree $T$ of order $n$ containing $P_4$, it holds that
$\prop_4(T) = n-1 $ if and only if $T$ is a subgraph of a multi-corona of $P_3$ whose
middle vertex has degree $2$.
\end{theorem}

\section{Coloring Trees Without a Rainbow $P_5$}  \label{s:p5}

In this section we determine the trees with the minimum value of the parameters for colorings without a rainbow $P_5$.

\begin{theorem}  \label{t:p5lower}
For a tree $T$ of order $n\ge 3$, it holds that $\upp_5(T) \ge \prop_5(T) \ge (n+3)/2$, and this is best possible.
\end{theorem}
\begin{Proof}
Since a coloring without a rainbow $P_4$ also doesn't have a rainbow $P_5$, we know from 
Theorem~\ref{t:first} that $\prop_5(T) \ge (n+2)/2$. Thus it suffices to show that 
achieving $(n+2)/2$ is not possible. 

As noted in the lead-in to Theorem~\ref{t:first}, 
one obtains a proper no-rainbow-$P_4$ coloring in $T$ by choosing one partite set $X$ and giving every vertex in  
$X$ its own unique color while giving all vertices in the other partite set $Y$ the same color.
So $\prop_5(T) \ge \prop_4(T) > n/2+1$ unless both partite sets have size $n/2$. So assume that is the case.

Choose some end-vertex $y$; say it is in $Y$. Then use the same coloring as above,
except that $y$ also gets its own unique color. Then the total number of colors is $n/2+2$.
Every $P_5$ that contains $y$ contains
two other vertices of $Y$, and hence remains not rainbow.
That is, the no-rainbow-$P_5$ number of $T$ is more than $(n+2)/2$.
\end{Proof}

\begin{figure}[h]
\centerline{\includegraphics{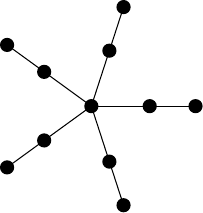}}
\caption{The octopus $O_5$}
\label{f:octopus}
\end{figure}

The value in Theorem~\ref{t:p5lower} is achieved by the octopus $O_b$ produced by taking the star with $b$ edges and subdividing every edge.
Figure~\ref{f:octopus} shows $O_5$.

\begin{lemma}
For the octopus $O_b$ with $b\ge 2$, it holds that $\upp_5(O_b) = \prop_5 (O_b) = b+2$.
\end{lemma}
\begin{Proof}
It is immediate that a minimum $P_5$-thwarting set is obtained by taking one edge from $b-1$ of the arms;
thus $\upp_5(O_b) = (2b+1) - (b-1) = b+2$. The value of $\prop_5(O_b)$ then follows by the above theorem.
\end{Proof}

We conclude this section by showing that the octopus is the unique extremal graph for both parameters:

\begin{theorem}
For odd $n\ge 5$, the octopus is the unique tree $T$ of order $n$ with $\prop_5(T) = (n+3)/2$ (and hence unique for $\upp_5$ too).
\end{theorem}
\begin{Proof}
It is immediate that the only tree of order $5$ with $\prop_5(T) < n$ is $P_5$ itself (which is the same as $O_2$). So
assume the tree $T$ has odd order $n\ge 7$ with bipartition $(X,Y)$ where $|X|>|Y|$.
It is immediate that $T$ is not a star.
As noted in the proof of Theorem~\ref{t:p5lower},  if there is an end-vertex in $Y$ then 
$\prop_5(T) \ge |X|+2$. It follows that all end-vertices of $T$ must be in $X$.
In particular, the diameter is at least four.

Consider a longest path in $T$, say starting $abcde$. 
Since all end-vertices are in $X$, no neighbor of $c$ is an end-vertex.
Let $A$ denote the set of
all vertices closer to $c$ than to $d$, other than $c$ itself.
Let $T' = T - A$. Note that $c$ is an end-vertex in tree $T'$. 
Any valid coloring of $T'$ can be extended to one of~$T$ by giving
each vertex of $A$ adjacent to $c$ the color of vertex $d$, and giving 
each vertex of $A$ nonadjacent to $c$ its own unique color.
Hence $\prop_5(T) \ge \prop_5(T') + |A|/2$.

It follows that $\prop_5(T) > (n+3)/2$, unless
$\prop_5(T') = (|T'|+3)/2$ 
and exactly half the vertices of $A$ are neighbors of $c$. 
Suppose that $T'$ has at least five vertices. Then, 
by the inductive hypothesis, the subtree $T'$ 
is an octopus with center $e$ while vertex $c$ is an end-vertex thereof.
But such a graph has $\prop_5 \ge (n+5)/2$: give every vertex of $X$ 
its own unique color except that $c$ and $e$ share colors,
give $d$ its own unique color, all other neighbors of $c$ share one color,
and all other neighbors of $e$ share another color. The number of 
colors used is $(|X|-1)+3$, a contradiction.
Hence in fact $T'$ has  
order $3$. Since $c$ is 
an end-vertex of $T'$, it follows that $T$ is an octopus.
\end{Proof}

\section{Conclusion}  \label{s:end}

For future work, one open question is to determine all trees with the minimum value of~$\prop_4$. 
It would also be of interest to consider bounds for other graph families, such as planar graphs or regular graphs;
for example, in~\cite{GXsecond} it is conjectured that $\prop_4 (G) \le n/2+1$ for every connected cubic graph $G$ of order~$n$.
And, of course, it would also be worthwhile to establish analogous bounds for $\upp_k$ and $\prop_k$ for larger $k$.

\end{document}